\documentclass[12pt, a4paper]{amsart}
\usepackage{amsmath}
\usepackage{geometry,amsthm,graphics,tabularx,amssymb,shapepar}
\usepackage{amscd}
\usepackage[all,2cell,dvips]{xy}

\newcommand{\Ad}{{\mathrm{Ad}}}

\newcommand{\End}{{\mathrm{End}}}

\newcommand{\GL}{{\mathrm{GL}}}

\newcommand{\rk}{{\mathrm{k}}}

\newcommand{\SO}{{\mathrm{SO}}}

\newcommand{\Sp}{{\mathrm{Sp}}}

\newcommand{\tr}{{\mathrm{tr}}}

\newcommand{\vsp}{{\vspace{0.2in}}}

\newcommand{\oO}{\operatorname{O}}
\newcommand{\oS}{\operatorname{S}}

\newcommand{\oU}{\operatorname{U}}

\renewcommand{\sl}{\mathfrak s \mathfrak l}

\renewcommand{\rk}{\mathrm k}

\newcommand{\C}{\mathbb{C}}
\newcommand{\R}{\mathbb R}

\newcommand{\ve}{{\vee}}

\newcommand{\la}{\langle}
\newcommand{\ra}{\rangle}

\newcommand{\be}{\begin {equation}}
\newcommand{\ee}{\end {equation}}
\newcommand{\bee}{\begin {equation*}}
\newcommand{\eee}{\end {equation*}}

\theoremstyle{Theorem}

\theoremstyle{Theorem}

\theoremstyle{Theorem}

\theoremstyle{Theorem}
\newtheorem{prp}{Proposition}[section]

\newtheorem{lemp}[prp]{Lemma}
\newtheorem{thmp}[prp]{Theorem}
\newtheorem{prpp}[prp]{Proposition}

\theoremstyle{Plain}

\theoremstyle{Definition}

\begin{document}

\title{Notes on MVW-extensions}

\author{Binyong Sun}

\address{Academy of Mathematics and Systems Science, Chinese Academy of
Sciences, Beijing, 100190, China}

\email{sun@math.ac.cn}




\begin{abstract}
We review certain basic geometric and analytic results concerning
MVW-extensions of classical groups, following
Moeglin-Vigneras-Waldspurger. The related results for Jacobi groups,
metaplectic groups, and special orthogonal groups are also included.
\end{abstract}

 \maketitle


\section{Preliminaries}\label{intro}

Let $\rk$ be a field, say, of characteristic zero. It is well known
that for every square matrix $x$ with coefficients in $\rk$, its
transpose $x^t$ is conjugate to $x$ (namely, there is an invertible
matrix $g$ with coefficients in $\rk$ such that $x^t=gxg^{-1}$). In
this note, we review some basic results on classical groups and
other related groups, which are closely related to this simple fact.

\subsection{$\epsilon$-Hermitian modules}\label{hermitianm}
We introduce the following terminologies and notations in order to
treat all classical groups uniformly. We find that this general
setting is very convenient when we apply Harish-Chandra descent in
the proof. Let $A$ be commutative semisimple finite-dimensional
$\rk$-algebra. It is thus a finite product of finite field
extensions of $\rk$. Let $\tau$ be a $\rk$-algebra involution on
$A$. We call such a pair $(A,\tau)$ a commutative involutive algebra
(over $\rk$). It is said to be simple if it is nonzero, and every
$\tau$-stable ideal of $A$ is either zero or $A$. Denote by $A^+$
the algebra of $\tau$-invariant elements in $A$. Then $(A,\tau)$ is
simple if and only if $A^+$ is a field. When this is the case,
either $A=A^+$, or $(A,\tau)$ is isomorphic to one of the
followings:
\[
  (A^+\times A^+, \tau_{A^+}), \quad(\textrm{a quadratic field
  extension of $A^+$, the nontrivial Galois element}),
\]
where $\tau_{A^+}$ is the coordinate exchanging map. Every
commutative involutive algebra is uniquely (up to ordering) a
product of simple ones.

Let $\epsilon=\pm 1$ and let $E$ be an $\epsilon$-Hermitian
$A$-module, namely it is a finitely generated $A$-module, equipped
with a non-degenerate $\rk$-bilinear map
\[
  \la\,,\,\ra_E:E\times E\rightarrow A
\]
satisfying
\[
     \la u,v\ra_E=\epsilon\la v,u\ra_E^\tau, \quad \la au,v\ra_E=a\la u,
     v\ra_E,\quad a\in A,\, u,v\in E.
\]

Denote by $\oU(E)$ the group of all $A$-module automorphisms of $E$
which preserve the form $\la\,,\,\ra_E$. When $(A,\tau)$ is simple,
$E$ is free as an $A$-module, and $\oU(E)$ is a classical group as
in the following table:
\begin{table}[h]
\caption{}
\centering 
\begin{tabular}{c c c c c c c} 
\hline
$A$ & \vline & $A^+$ & $A^+\times A^+$ & quadratic filed extension\\
\hline
$\epsilon=1$ & \vline & orthogonal group & general linear group  & unitary group\\ 
\hline 
$\epsilon=-1$ & \vline & symplectic group & general linear group  & unitary group \\
\hline 
\end{tabular}
\label{table:nonlin} 
\end{table}

In general, write
\begin{equation}\label{decoma}
  (A,\tau)=(A_1,\tau_1)\times (A_2,\tau_2)\times\cdots \times (A_r,\tau_r)
\end{equation}
as a product of simple commutative involutive algebras. Then
$E_i:=A_i\otimes_A E$ is obviously an $\epsilon$-Hermitian
$A_i$-module. We have that
\begin{equation}\label{decome}
  E=E_1\times E_2\times \cdots \times E_r,
\end{equation}
and
\[
  \oU(E)=\oU(E_1)\times \oU(E_2)\times \cdots \oU(E_r).
\]

We say that $E$ is simple if it is nonzero, and every non-degenerate
$A$-submodule of it is either zero or $E$. Every
$\epsilon$-Hermitian $A$-module is isomorphic to an orthogonal sum
of simple ones.

For every
\[
 a\in (A^{\times})^{\tau=\epsilon}:=\{a\in A^\times\mid a^\tau=\epsilon
 a\},
\]
write $A(a):=A$ as an $A$-module, equipped with the form
\[
  \la u,v\ra_{A(a)}=auv^\tau\in A,\quad u,v\in A(a).
\]
Then $A(a)$ is an $\epsilon$-Hermitian $A$-module,
and $A(a)$ is isomorphic to $A(a')$ if and only if
\[
 a'a^{-1}=bb^\tau \,\,\textrm{  for some }b\in A^\times.
\]

The following classification of simple $\epsilon$-Hermitian
$A$-modules is obvious.

\begin{prp}\label{simmpleh}
Assume that $(A,\tau)$ is simple.

(a) If $A=A^+\times A^+$, then there is a unique simple
$\epsilon$-Hermitian $A$-module (up to isomorphism). It has rank
$1$.

(b) If $A=A^+$ and $\epsilon$=-1, then there is a unique simple
$\epsilon$-Hermitian $A$-module. It has rank $2$.

(c) If $A=A^+$ and $\epsilon=1$, or $A$ is a quadratic field
extension of $A^+$, then every simple $\epsilon$-Hermitian
$A$-module is of the form $A(a)$ ($a\in
(A^{\times})^{\tau=\epsilon}$).
\end{prp}

Write $E_\tau:=E$ as a $\rk$-vector space, and for every $v\in E$,
write $v_\tau:=v$, viewed as a vector in $E_\tau$. Then $E_\tau$ is
an $\epsilon$-Hermitian $A$-module under the scalar multiplication
\[
  a v_\tau:=(a^\tau v)_\tau, \quad a\in A, v\in E,
\]
and the form
\[
  \la u_\tau, v_\tau\ra_{E_\tau}:=\la v,u\ra_E, \quad u,v \in E.
\]

\begin{prpp}\label{irreh}
The $\epsilon$-Hermitian $A$-modules $E_\tau$ and $E$ are isomorphic
to each other.
\end{prpp}
\begin{proof}
Without loss of generality, assume that $(A,\tau)$ is simple and $E$
is simple. Then the proposition follows trivially from Proposition
\ref{simmpleh}.
\end{proof}

\subsection{Harish-Chandra descent}
Let $E$ be an $\epsilon$-Hermitian $A$-modules. Define an involution
on $\End_A(E)$, which is still denoted by $\tau$,  by requiring that
\[
  \la xu,v\ra_E=\la u, x^\tau v\ra_E, \quad x\in \End_A(E), \,\,u,v\in E.
\]

Let $s$ be a semisimple element of $\End_A(E)$ (that is, it is
semisimple as a $\rk$-linear operator). Assume it is normal in the
sense that $s^\tau s=s s^\tau$. Denote by $A_s$ the subalgebra of
$\End_A(E)$ generated by $s$, $s^\tau$ and scalar multiplications by
$A$. It is $\tau$-stable and $(A_s,\tau)$ is a commutative
involutive algebra. Write $E_s:=E$, viewed as an $A_s$-module.
Define a $\rk$-bilinear map
\[
   \la\,,\,\ra_{E_s}: E_s \times E_s \rightarrow A_s
\]
by requiring that
\[
    \tr_{A_s/\rk}(a\la u, v\ra_{E_s})= \tr_{A/\rk}(\la a u, v\ra_{E}), \quad u,v\in E, \,\,a\in A_s.
\]
Then $E_s$ becomes an $\epsilon$-Hermitian $A_s$-modules. When $s\in
\oU(E)$, geometric and analytic problems on $\oU(E)$ around $s$ are
often reduced to that on $\oU(E_s)$. The procedure is called
Harish-Chandra descent.

\subsection{$\epsilon$-Hermitian $\sl_2$-modules}

Let $A$, $\tau$, $A^+$, $\epsilon$ be as before. Let $E$ be an
$\epsilon$-Hermitian $(\sl_2,A)$-module, namely it is an
$\epsilon$-Hermitian $A$-module, equipped with a Lie algebra action
\[
  \sl_2(\rk)\times E\rightarrow E, \quad x,v\mapsto xv,
\]
which is $\rk$-linear on the first factor,  $A$-linear on the second
factor, and satisfies
\[
  \la xu,v\ra_E+\la u,xv\ra_E=0, \quad x\in \sl_2(\rk),\,u,v\in E.
\]
We say that $E$ is irreducible if it is nonzero, and every
$\sl_2(\rk)$-stable  $A$-submodule of it is either zero or $E$. We
say that $E$ is simple if it is nonzero, and every
$\sl_2(\rk)$-stable non-degenerate $A$-submodule of it is either zero
or $E$. Every $\epsilon$-Hermitian $(\sl_2, A)$-module is isomorphic
to an orthogonal sum of simple ones.

Write
\[
  \mathbf h:= \left[
                   \begin{array}{cc} 1&0\\ 0&-1\\
                   \end{array}
  \right], \quad
  \mathbf e:= \left[
                   \begin{array}{cc} 0&1\\ 0&0\\
                   \end{array}
  \right],\quad
  \mathbf f:=\left[
                   \begin{array}{cc} 0&0\\ 1&0\\
                   \end{array}
  \right],
\]
which form a basis of the Lie algebra $\sl_2(\rk)$. The following
lemma is also obvious.

\begin{lemp}
Assume that $(A,\tau)$ is simple. If $A=A^+$ and $\epsilon=1$, or
$A$ is a quadratic field extension of $A^+$, then for every positive
odd integer $2d-1$ and every $a\in (A^{\times})^{\tau=\epsilon}$,
there is a unique (up to isomorphism) simple $\epsilon$-Hermitian
$(\sl_2,A)$-module with the following property: it has rank $2d-1$
and its $\mathbf h$-invariant vectors form an $\epsilon$-Hermitian
$A$-module which is isomorphic to $A(a)$.
\end{lemp}

Write $A^{2d-1}(a)$ for this simple $\epsilon$-Hermitian
$(\sl_2,A)$-module. Then $A^{2d-1}(a)$ is isomorphic to
$A^{2d-1}(a')$ if and only if
\[
 a'a^{-1}=bb^\tau \,\,\textrm{  for some }b\in A^\times.
\]

The following proposition, which follows easily from the
representation theory of $\sl_2$, classifies simple
$\epsilon$-Hermitian $(\sl_2,A)$-modules.

\begin{prpp}\label{simmplesl}Assume that $(A,\tau)$ is simple.

(a)Further assume that $A=A^+\times A^+$. Then for every positive
integer $d$, there is a unique simple $\epsilon$-Hermitian
$(\sl_2,A)$-module of rank $d$.

(b)Further assume that $A$ is a quadratic field extension of $A^+$.
Then for every positive even integer $2d$, there is a unique simple
$\epsilon$-Hermitian $(\sl_2,A)$-module of rank $2d$; every simple
$\epsilon$-Hermitian $(\sl_2,A)$-module of odd rank is of the form
$A^{2d-1}(a)$ ($a\in (A^{\times})^{\tau=\epsilon}$).

(c)Further assume that $A=A^+$ and $\epsilon$=-1. Then for every
positive even integer $d$, there is a unique simple
$\epsilon$-Hermitian $(\sl_2,A)$-module of rank $2d$; for every
positive odd integer $d$, there are exactly two simple
$\epsilon$-Hermitian $(\sl_2,A)$-module of rank $2d$, one is
irreducible and the other is reducible.

(d)Further assume that $A=A^+$ and $\epsilon$=1. Then for every
positive integer $d$, there is a unique simple $\epsilon$-Hermitian
$(\sl_2,A)$-module of rank $4d$, and there is no  simple
$\epsilon$-Hermitian $(\sl_2,A)$-module of rank $4d-2$; every simple
$\epsilon$-Hermitian $(\sl_2,A)$-module of odd rank is of the form
$A^{2d-1}(a)$ ($a\in A^{\times}$).

\end{prpp}

Recall that $E$ is an $\epsilon$-Hermitian $(\sl_2,A)$-module. We
define an $\epsilon$-Hermitian $(\sl_2,A)$-module $E_\tau$ as
follows. As an $\epsilon$-Hermitian $A$-module, it is as in Section
\ref{hermitianm}. The $\sl_2(\rk)$-action is given by
\[
  \mathbf h v_\tau:=(\mathbf h v)_\tau, \quad \mathbf e
  v_\tau:=-(\mathbf e v)_\tau, \quad \mathbf f v_\tau:=-(\mathbf f
  v)_\tau, \quad v\in E.
\]

\begin{prpp}
The $\epsilon$-Hermitian $(\sl_2,A)$-modules $E_\tau$ and $E$ are
isomorphic to each other.
\end{prpp}
\begin{proof}
This follows trivially from Proposition \ref{simmplesl}.
\end{proof}

\section{Geometric results}

\subsection{Classical groups}

Following Moeglin-Vigneras-Waldspurger (\cite[Proposition
4.I.2]{MVW87}), we extend $\oU(E)$ to a larger group, which is
denoted by $\breve{\oU}(E)$, and is defined to be the subgroup of
$\GL(E_\rk)\times \{\pm 1\}$ consisting of pairs $(g,\delta)$ such
that either
\[
  \delta=1 \quad\textrm{and } \la gu,gv\ra_E=\la u,v\ra_E,\quad  u,v\in  E,
\]
or
\[
  \delta=-1 \quad\textrm{and }\la gu,gv\ra_E=\la v,u\ra_E,\quad  u,v\in E.
\]
Here $E_\rk:=E$, viewed as a $\rk$-vector space. Every $g\in
\GL(E_\rk)$ is automatically $A$-linear if $(g,1)\in
\breve{\oU}(E)$, and is conjugate $A$-linear (with respect to
$\tau$) if $(g,-1)\in \breve{\oU}(E)$. We call $\breve{\oU}(E)$ the
MVW-extension of $\oU(E)$. Proposition \ref{irreh} amounts to saying
that the projection map $\breve{\oU}(E)\rightarrow \{\pm 1\}$ is
surjective. Therefore, we have a short exact sequence
\[
     1\rightarrow\oU(E)\rightarrow
   \breve{\oU}(E)\rightarrow
    \{\pm 1\}\rightarrow 1.
\]

The following basic fact of classical groups is a part of
\cite[Proposition 4.I.2]{MVW87}. With the preparation of Section 1,
we sketch a short proof here.

\begin{thmp}\label{geo} For every $x\in \oU(E)$, there is an element $\breve
g\in \breve{\oU}(E)\setminus \oU(E)$ such that
$\breve{g}x\breve{g}^{-1}=x^{-1}$.
\end{thmp}

\begin{proof}
By using Jordan decomposition and Harish-Chandra descent, we may
(and do) assume that $x$ is unipotent.  By Jacobson-Morozov Theorem,
we choose an action of $\sl_2(\rk)$ on $E$ such that it makes $E$ an
$\epsilon$-Hermitian $(\sl_2,A)$-module, and that the exponential of
the action of $\mathbf e$ coincides with $x$. Then the theorem
follows from Proposition \ref{simmplesl}.
\end{proof}

Moeglin-Vigneras-Waldspurger also prove the Lie algebra analog of
Theorem \ref{geo}, namely, for every $x$ in the Lie algebra of
$\oU(E)$, there is an element $\breve g\in \breve{\oU}(E)\setminus
\oU(E)$ such that $\Ad_{\breve{g}}x=-x$. When $\oU(E)$ is a general
linear group, this is just the simple fact mentioned at the
beginning of this note.

Similar to Theorem \ref{geo}, the following Theorems
\ref{geo2}-\ref{geo6} (as well as \cite[Proposition 4.I.2]{MVW87})
can be proved by using Harish-Chanda descent and representation
theory of $\sl_2$.  We leave the details to the interested reader.

\subsection{Jacobi groups}\label{sjacobi}
In the remaining part of this note, we assume for simplicity that $(A,\tau)$ is
simple. Denote by $L$ a free $A$-submodule of $E$ of
rank one (if $E$ is nonzero), and by $L^+$ a cyclic $A^+$-subspace
of $L$ which generates $L$ as an $A$-module.

Write $\oU_L(E)$ for the subgroup of $\oU(E)$ fixing $L$ point-wise,
and write
\[
  \breve{\oU}_{L^+}(E):=\{(g,\delta)\in\breve{\oU}(E)\mid  g
  \textrm{ fixes $L^+$ point-wise}\}.
\]
It contains $\oU_L(E)$ as a subgroup of index two.

Similar to Theorem \ref{geo}, we have

\begin{thmp}\label{geo2} Assume that $L$ is totally isotropic. Then for every $x\in \oU_L(E)$, there is an element
$\breve g\in \breve{\oU}_{L^+}(E)\setminus \oU_L(E)$ such that
$\breve{g}x\breve{g}^{-1}=x^{-1}$.
\end{thmp}

It seems that Theorem \ref{geo2} (and Theorem {\ref{geo6} of
Section \ref{smet}) are new. If $L$ is not totally isotropic, then it
is non-degenerate and Theorem \ref{geo2} also holds (which is a
restatement of Theorem \ref{geo}).

\subsection{Special orthogonal groups}

Assume that $A=A^+$ and $\epsilon=1$. Following Waldspurger (\cite{Wald09}), we define
\[
 \oS\!\breve{\oO}(E):=\left\{ (g,\delta)\in \oO(E)\times \{\pm 1\}\mid \det(g)=\delta^{\left[\frac{\dim E +1}{2} \right]}\right \}.
\]
It contains the special orthogonal group $\SO(E)$ as a subgroup of
index two. By convention, $\oS\!\breve{\oO}(E):=\oO(E)\times \{\pm
1\}=\{\pm 1\}$ if $E=\{0\}$.

\begin{thmp}\label{geo3} (cf. \cite{Wald09}) For every $x\in \SO(E)$, there is an element
$\breve g\in \oS\!\breve{\oO}(E)\setminus \SO(E)$ such that
$\breve{g}x\breve{g}^{-1}=x^{-1}$.
\end{thmp}

\subsection{Metaplectic groups}\label{smet}

In the remaining part  of this note, we assume that $\rk$ is a local field of
characteristic zero. In this subsection further assume that $\epsilon=-1$. Write $E_\rk:=E$,
viewed as a $\rk$-symplectic space under the form
\[
  \la u,v\ra_{E_\rk}:=\tr_{A/\rk}(\la u,v\ra_E).
\]
Denote by
\begin{equation}\label{meta}
   1\rightarrow \{\pm 1\}\rightarrow
   \widetilde{\Sp}(E_\rk)\rightarrow
   \Sp(E_\rk)\rightarrow 1
\end{equation}
the metaplectic cover of the symplectic group $\Sp(E_\rk)$. This is
a (topologically) exact sequence of locally compact topological
groups. It splits when either $E=0$ or $\rk=\C$. Otherwise, this is
the unique non-split (topological) central extension of $\Sp(E_\rk)$
by $\{\pm 1\}$ (cf. \cite[Theorem 10.4]{Moo68}).

Note that $\breve{\Sp}(E_\rk):=\breve{\oU}(E_\rk)$ equals to the
subgroup of $\operatorname{GSp}(E_\rk)$ with similitudes $\pm 1$. It
is shown  in \cite[Page 36]{MVW87} that there is a unique action
\begin{equation}\label{liftact0}
   \Ad:  \breve{\Sp}(E_\rk)\times \widetilde{\Sp}(E_\rk)\rightarrow \widetilde{\Sp}(E_\rk)
\end{equation}
of $\breve{\Sp}(E_\rk)$ as group automorphisms on
$\widetilde{\Sp}(E_\rk)$ which lifts the adjoint action
\[
  \Ad: \breve{\Sp}(E_\rk)\times\Sp(E_\rk)\rightarrow \Sp(E_\rk)
\]
and fixes the central element $-1\in \widetilde{\Sp}(E_\rk)$.

Denote by $\widetilde{\oU}(E)$ the double cover of $\oU(E)\subset
\Sp(E_\rk)$ induced by the cover (\ref{meta}). Then the action
(\ref{liftact0}) restricts to an action
\[
 \Ad:
\breve{\oU}(E)\times \widetilde{\oU}(E)\rightarrow
\widetilde{\oU}(E).
\]

\begin{thmp}\label{contras}
For every $x\in \widetilde \oU(E)$, there is an element $\breve g\in
\breve{\oU}(E)\setminus \oU(E)$ such that $\Ad_{\breve g} x=x^{-1}$.
\end{thmp}

The most interesting case is when $\oU(E)$ is a symplectic group. Then Theorem \ref{contras} is proved for semisimple elements
in \cite[Proposition 4.I.8]{MVW87} and for general elements in
\cite[Theorem 1.1]{FS}.

\vsp

Recall $L$ and $L^+$ from Section \ref{sjacobi}. Denote by $\widetilde \oU_L(E)$ the double cover of $\oU_L(E)\subset
\Sp(E_\rk)$ induced by the cover (\ref{meta}). The action
(\ref{liftact0}) restricts to an action
\[
 \Ad:
\breve{\oU}_{L^+}(E)\times \widetilde{\oU}_L(E)\rightarrow
\widetilde{\oU}_L(E).
\]

\begin{thmp}\label{geo6} Assume that $L$ is totally isotropic. Then for every $x\in \widetilde{\oU}_L(E)$, there is an element
$\breve g\in \breve{\oU}_{L^+}(E)\setminus \oU_L(E)$ such that
$\Ad_{\breve g} x=x^{-1}$.
\end{thmp}

\section{Analytic results}

Recall that $\rk$ is a local field of characteristic
zero, $(A,\tau)$ is simple, and $E$ is an $\epsilon$-Hermitian
$A$-module.

\subsection{Contragredient representations}

Let $G$ denote one of the following groups:
\[
  \oU(E), \quad  \SO(E)\,\,(\textrm{when $A=A^+$ and $\epsilon=-1$}), \quad \widetilde{\oU}(E)\,\, (\textrm{when }\epsilon=-1),
\]
or one of the following groups when there is a totally isotropic
rank one free $A$-submodule $L$ of $E$:
\[
   \oU_L(E), \quad\widetilde{\oU}_L(E)\,\, (\textrm{when }\epsilon=-1).
\]
Let $\breve g$ be respectively an element of
\[
  \breve{\oU}(E)\setminus \oU(E), \quad \oS\!\breve{\oO}(E)\setminus \SO(E), \quad\breve{\oU}(E)\setminus
  \oU(E),
\]
or
\[
  \breve{\oU}_{L^+}(E)\setminus \oU_L(E),  \quad\breve{\oU}_{L^+}(E)\setminus
  \oU_L(E).
\]
Here and as before, $L^+$ is a one-dimensional $A^+$-subspace which
generates $L$ as an $A$-module. In all cases, we have a group
automorphism $\Ad_{\breve g}:G\rightarrow G$.

\begin{thmp}\label{contraf}
For every  invariant (under the adjoint action of $G$) generalized
function $f$ on $G$, one has that  $f(\Ad_{\breve{g}}x)=f(x^{-1})$
(as generalized functions on $G$).
\end{thmp}

For the usual notion of generalized functions, see \cite[Section
2]{Sun08} (non-archimedean case) and \cite[Section 2.1]{JSZ}
(archimedean case), for example.  By the localization principle of
Bernstein and Zelevinsky (\cite[Theorem 6.9]{BZ76}), Theorem
\ref{geo}-\ref{geo6} implies Theorem \ref{contraf} in the
non-archimedean case. In both archimedean and non-archimedean
cases, Theorem \ref{contraf} is implied by Theorems \ref{strongi}
and \ref{strongii} of the next two subsections.

\vsp

When $\rk$ is non-archimedean, Theorem \ref{contraf} implies that
for every irreducible admissible smooth representation $\pi$ of $G$,
its contragredient $\pi^\ve$ is isomorphic to its twist $x\mapsto
\pi(\Ad_{\breve g}x)$. Certain archimedean analog of this fact also
holds. But it is less satisfactory due to the lack of a suitable
notion of ``admissible representations" for non-reductive real
groups.

\subsection{Multiplicity one theorem I}\label{seci}

Only in this subsection, assume that $L$ is non-degenerate. (This is
imposable in the symplectic case.) Denote by $E_0$ the orthogonal complement of $L$ in $E$.

Let $G$ denote one of the following groups:
\[
  \oU(E), \quad  \SO(E)\,\,(\textrm{when $A=A^+$ and $\epsilon=-1$}), \quad\widetilde{\oU}(E)\,\, (\textrm{when }\epsilon=-1),
\]
and let $G_L$ denote its respective subgroup
\[
  \oU(E_0), \quad  \SO(E_0), \quad\widetilde{\oU}(E_0).
\]
When $G=\SO(E)$, let $\breve g$ be an element
\[
  \oS\!\breve{\oO}(E_0)\setminus \SO(E_0)\subset \oO(E_0)\times \{\pm 1\}\subset \oO(E)\times \{\pm 1\}\supset G,
\]
and in all other cases, let $\breve g$ be an element of
of
\[
   \breve{\oU}(E_0)\setminus \oU(E_0)=\breve{\oU}_{L^+}(E)\setminus \oU_{L}(E)\subset \breve \oU(E).
\]
As before, we have a group automorphism $\Ad_{\breve g}:G\rightarrow G$.

\begin{thmp}\label{strongi}
For every generalized function $f$ on $G$ which is invariant under
the adjoint action of $G_L$, one has that
$f(\Ad_{\breve{g}}x)=f(x^{-1})$.
\end{thmp}

In the non-archimedean case, Theorem \ref{strongi} is proved
by Aizenbud-Gourevitch-Rallis-Schiffmann in \cite{AGRS} (except for the case of special orthogonal groups, which is proved by Waldspurger in \cite{Wald09}). In the
archimedean case, it is proved by Sun-Zhu in \cite{SZ} (and
independently by Aizenbud-Gourevitch in \cite{AG} for general linear groups).

By Gelfand-Kazhdan criteria, Theorem \ref{strongi} implies that
$(G,G_L)$ is a ``multiplicity one pair" (see \cite{GGP, AGRS, AG,
SZ}). This multiplicity one theorem has been expected by Bernstein
and Rallis since 1980's.

\subsection{Multiplicity one theorem II}
Assume that $L$ is totally isotropic in this subsection. Let $L'$ be
another  totally isotropic rank one free $A$-submodule of $E$ which
is dual to $L$ under the form $\la\,,\,\ra_E$. Then
\begin{equation}\label{decome}
  E:=L\oplus E_0\oplus L',
\end{equation}
where $E_0$ is the orthogonal complement of $L\oplus L'$ in $E$.

Let $G_L$ denote one of the following groups
\[
   \oU_L(E), \quad\widetilde{\oU}_L(E)\,\, (\textrm{when }\epsilon=-1).
\]
Let $G_0$  denote its respective subgroup $\oU(E_0)$ or $\widetilde{\oU}(E_0)$. Let $\breve g=(g,-1)\in \breve{\oU}(E_0)\setminus \oU(E_0)$. View it as an element of $\breve{\oU}(E)\setminus \oU(E)$ by extending $g$ to a $\tau$-conjugate linear automorphism of $E$, preserving the decomposition \eqref{decome} and fixing $L^+$ point-wise.  Again, we have a group automorphism $\Ad_{\breve g}:G_L\rightarrow
G_L$.

\begin{thmp}\label{strongii}
For every generalized function $f$ on $G_L$ which is invariant under
the adjoint action of $G_0$, one has that
$f(\Ad_{\breve{g}} x)=f(x^{-1})$.
\end{thmp}

For a proof of Theorem \ref{strongii}, see \cite{Sun08, Di, SZ}.
Similar to Theorem \ref{strongi}, Theorem \ref{strongii}
implies that the pair $(G_L, G_0)$ is a ``multiplicity one pair". This multiplicity one theorem was expected by Prasad, at least for symplectic groups (\cite[Page 20]{Pr96}). 

\vsp

\noindent {\bf Remarks}:  In fact, the metaplectic cases of Theorem \ref{strongi} and Theorem \ref{strongii} were not treated in the literature. However, the available method works as well.

\subsection{The group $\oS\!\breve{\oO}(E)\ltimes E$}
Assume that $A=A^+$  and  $\epsilon=1$. Let $\oS\!\breve{\oO}(E)$ act on $E$ by $(g,\delta)v:=gv$, and we
form the semidirect product $\oS\!\breve{\oO}(E)\ltimes E$. In general, the
desired geometric property does not hold for this group. For example, if $E$ is split and has dimension 2, and $x\in E\subset \SO(E)\ltimes E$ is a nonzero isotropic vector, then there is no
$\breve g\in (\oS\!\breve{\oO}(E)\ltimes E) \setminus (\SO(E)\ltimes
E)$ such that $\breve gx \breve g^{-1}=x^{-1}$. However, the corresponding analytic result still holds:
\begin{thmp}\label{contraso}
For every invariant (under the adjoint action of
$\SO(E)\ltimes E$) generalized function $f$ on $\SO(E)\ltimes E$, and for every element of $\breve g$ of $(\oS\!\breve{\oO}(E)\ltimes E) \setminus (\SO(E)\ltimes E)$, one has that
$f(\breve{g}x\breve{g}^{-1})=f(x^{-1})$.
\end{thmp}

This is much weaker than the following
\begin{thmp}\label{moso2}
Let $\breve g\in \oS\!\breve{\oO}(E)\setminus \SO(E)$.  Then for
every generalized function $f$ on $\SO(E)\ltimes E$ which is
invariant under the adjoint action of $\SO(E)$, one has that
$f(\breve gx\breve g^{-1})=f(x^{-1})$.
\end{thmp}

Theorem \ref{moso2} can be proved by using the methods and results of Aizenbud-Gourevitch-Rallis-Schiffmann and Sun-Zhu (cf. \cite{Di, Sun08, Wald09, SZ}).

\vsp \noindent Acknowledgements: the author thanks Chee Whye Chin,
Wee Teck Gan, Hung Yean Loke and Chen-Bo Zhu for helpful
discussions. The work was supported by NSFC grants 10801126 and
10931006.

\end{document}